\documentclass[a4paper,11pt]{article}
\usepackage{amssymb,amsmath, amsthm}
\usepackage[all]{xy}

\title
{Explicit Kummer surface formulas for arbitrary characteristic}

\author{Jan Steffen M\"uller 
 \thanks{Supported by DFG-grant STO 299/5-1}}

\newtheorem{thm}{Theorem}[section]
\newtheorem{prop}[thm]{Proposition}
\newtheorem{lemma}[thm]{Lemma}

\newcommand\Char{\mathop{\rm char}\nolimits}

\newcommand\isom{\cong}

\newcommand{\To}{\longrightarrow}
\newcommand{\BP}{{\mathbb P}}

\newcommand\into{\hookrightarrow}

\begin{document}
\maketitle
\begin{abstract}
If $C$ is a curve of genus~2 defined over a field $k$ and $J$ is its Jacobian, then we can associate a hypersurface
$K$ in $\mathbb{P}^3$ to $J$, called the Kummer surface of $J$. Flynn has made this construction 
explicit in the case that the characteristic of $k$ is not~2 and $C$ is given by a simplified 
equation. He has also given explicit versions of several maps defined on the Kummer surface and
shown how to perform arithmetic on $J$ using these maps. In this paper we generalize these results
to the case of arbitrary characteristic. 
\end{abstract}

\section{Introduction}\strut\\[1mm] \label{Intro}
If $C$ is a curve of genus~2 defined over a field $k$ and $J$ is its Jacobian variety, then $J$ can be regarded as a smooth projective variety embedded into $\BP^{15}$. However, in order to perform explicit arithmetic on $J$, this construction is not suitable, since computations in $\BP^{15}$ are too cumbersome. To remedy this, one can associate a hypersurface $K$ in $\BP^3$ to $J$, called the Kummer surface of $J$. It is the quotient of the Jacobian by the negation map and, although not an abelian variety itself, remnants of the group law on $J$ can be exhibited on $K$ and arithmetic on $J$ can be performed using its Kummer surface. The map from $J$ to $K$ corresponds to the map from an elliptic curve to $\BP^1$ that assigns to an affine point its $x$-coordinate.

In \cite{CasselsFlynn} Cassels and Flynn construct the Kummer surface associated to the Jacobian of a genus~2 curve $C$ defined over a field $k$ and several maps on it only in the special case 
\begin{equation}\label{SpecEq}
C: y^2=f(x)
\end{equation}
and leave the general case
\begin{equation}\label{GenEq}
 C:y^2+h(x)y=f(x),
\end{equation}
where $\deg(f)\le6$ and $\deg(h)\le 3$, as ``optional exercises for the reader'' (p.1). If $\Char(k)\ne 2$, then we can find a defining equation as in \eqref{SpecEq} for any genus~2 curve over $k$, but this is not true in the case $\Char(k)=2$. In the present paper these exercises will be tackled; the aim is to find expressions that work whenever the curve is given by a general equation as in \eqref{GenEq} over any field. For the special case $\Char(k)=2$ and $\deg(h)=2$ (which one can always reduce to if $\Char(k)=2$) such expressions have also been obtained independently by Duquesne (cf. \cite{SD}). However, they do not work over fields of different characteristic. Nevertheless, all expressions found in the present work specialize to those obtained in \cite{SD} in the case investigated there. It has to be noted that Duquesne also explicitly constructs the Jacobian, which has not been attempted by the author in the general setting.

Flynn presents an explicit theory of Kummer surfaces in the case that $\Char(k)\ne 2$ and $C$ is given by an equation as in \eqref{SpecEq}. In \cite{Flynnb} he introduces an explicit embedding of $K$ into $\BP^3$ and a defining equation for $K$. He also computes explicit expressions for a map $\delta$, corresponding to duplication on $J$, and two matrices $B$ (with entries certain biquadratic forms) and $W$, which corresponds to translation by a point of order~2, that help to use the Kummer surface to facilitate efficient arithmetic on the Jacobian. The discussion in \cite{CasselsFlynn} is essentially an exposition of his results.

It appears that so far no attempt has been made to construct such an explicit theory in the more general case that $k$ is any field and $C$ is a curve of genus~2 given by an equation as in \eqref{GenEq}. However, from general theory, we know that all relevant objects constructed by Flynn must have counterparts in this more general situation. 

Concerning possible applications, the author's original motivation was to improve the algorithm for the computation of canonical heights on Jacobians of curves of genus~2 defined over number fields or function fields introduced by Flynn and Smart in \cite{FlynnSmart} and modified by Stoll in \cite{StollH2}. Here one computes local heights on the Kummer coordinates of a point on the Jacobian for each valuation of the base field of the curve. The results of this work should help in two ways: On the one hand, a curve may have a model as in \eqref{GenEq} with much smaller coefficients than those of its simplified models \eqref{SpecEq} (this is similar to the elliptic curve situation). But it can be shown that the change in the local height caused by a change of model is given by a simple formula. Thus we may achieve a significant speed improvement, although it has to be taken into account that the formulas for duplication that the algorithm relies on might be more complicated. On the other hand, the simplified models are often not minimal for residue characteristic~2 and always have bad reduction there. Accordingly, the largest portion of the running time of the algorithm is usually spent on computing the local height for residue characteristic~2.

Another application is in the field of cryptography. Here one can use, more generally, the Jacobian of genus~2 curves defined over finite fields whose order is either a large prime or a large power of $2$. In the former case the Kummer surface can be used to speed up arithmetic on such Jacobians and hence to make cryptosystems based on them more practical; see for instance the discussion in \cite{SD2} and \cite{Gaudry1}. In the latter case, an explicit theory of the Kummer surface - as presented in this paper - can also be used, see \cite{SD3}. Recently, Gaudry and Lubicz (\cite{GaudryLubicz}) have also introduced formulas for the arithmetic of Kummer surfaces in characteristic~2 based on the theory of algebraic theta functions that are different from the results in the paper at hand and from those in \cite{SD}. They also show that such formulas can indeed be useful for cryptographic purposes. 

In this paper we give explicit expressions for the abovementioned objects that reduce to the expressions given in \cite{Flynnb} and \cite{SD} whenever those are valid. First, we briefly recall Flynn's construction in section \ref{classical}, before we give  an explicit embedding and a defining equation for the Kummer surface in general characteristic in section \ref{GenKum}. We do not attempt to copy Flynn's or Duquesne's approach (except for the computation of the matrix $W$ in characteristic~2 in section \ref{Wg}, see below), but rather make use of the fact that if $\Char(k)\ne 2$, then $K$ is isomorphic to the Kummer surface associated with the Jacobian of a genus~2 curve given by an equation as in \eqref{SpecEq}. We find this isomorphism explicitly and use it to map our formulas to $K$. Next, we slightly modify them and prove that these modified versions remain valid in characteristic~2. The duplication map is discussed in section \ref{Duplication} and the matrix $B$ of biquadratic forms is computed in section \ref{Biquads}. In fact the proofs are only sketched, but full details are given later in sections \ref{proofdelta} and \ref{proofB}. For the matrix $W$ we have to make a case distinction because the approach introduced above, although successful when $\Char(k)\ne 2$, cannot be easily modified so that it also works for fields of characteristic~2. Instead we compute $W$ directly in that case, using the same approach employed by Duquesne in the case $\deg(h)=2$.

Since the expressions are too complicated, we do not actually present all of them in this paper, but rather discuss how to obtain them. For the computations we have used a combination of the computer algebra systems MAGMA (\cite{Magma}) and Maple (\cite{Maple}).
\vspace{4mm}

\emph{Acknowledgements.}
\\I would like to thank my advisor Michael Stoll for the motivation for this work and for many valuable suggestions, in particular on exploiting the isomorphism in sections~4,5 and~6. I would also like to thank Victor Flynn for helpful discussions and Sylvain Duquesne for sending me a copy of his preprint \cite{SD}.

\section{The classical case}\strut\\[1mm] \label{classical}
In \cite{Flynnb} Flynn explicitly constructs the Kummer surface $K$ associated with the Jacobian of a genus~2 curve $C$ given by a model of the form
\[
C: y^2=f(x).
\]
Our intention is to generalize this construction to models of the form
\begin{equation}
 C:y^2+h(x)y=f(x).
\end{equation}
 We begin by reviewing Flynn's construction from \cite{Flynnb}. Let
\[
f(x)=f_0+f_1x+f_2x^2+f_3x^3+f_4x^4+f_5x^5+f_6x^6
\] 
be a polynomial in $k[x]$ without multiple factors, where $k$ is a field of characteristic $\Char(k)\ne 2$ and $\deg(f) = 6$. Then the affine equation 
\[
 y^2=f(x)
\]
defines a hyperelliptic curve $C$ of genus~2 defined over $k$. We denote its Jacobian variety by $J$. If we form the quotient of $J$ by the negation map, then we get another variety $K$, the Kummer surface associated with $J$. In \cite{Flynna} and \cite{Flynnb} Flynn explicitly realizes these objects. The Jacobian lives in $\BP^{15}$, whereas the Kummer surface lives in $\BP^{3}$, so explicit calculations are much more efficient on the Kummer surface. Since remnants of the group structure are preserved when one passes to the Kummer surface, these remnants can be used to obtain a feasible way of performing arithmetic on $J$.

The explicit embeddings of both the Jacobian and the Kummer surface can be found using a modified version of the classical theta-divisor on the Jacobian. The classical theta-divisor $\Theta$ over an algebraically closed field $k$ is defined to be the divisor on $J$ given by the image of $C$ under the embedding
\begin{eqnarray*}
 \iota : C & \into& J \\
   P_1 &\mapsto& [P_1-\infty],\\
\end{eqnarray*}
where the assumption that $k$ is algebraically closed means that we may assume $f_6=0$, and so $\infty$ is the unique point at infinity. On the other hand, if $k$ is not algebraically closed, then in the case $f_6\ne 0$ there are two branches $\infty^+$ and $\infty^-$ over the singular point at infinity. In this situation, we define $\Theta^+$ and $\Theta^-$ to be the images of $C$ under the embeddings
\begin{eqnarray*}
 \iota^+ : C & \into& J \\
 P_1 &\mapsto& [P_1-\infty^+]\\
\end{eqnarray*}
and
\begin{eqnarray*}
 \iota^- : C & \into& J \\
 P_1 &\mapsto& [P_1-\infty^-],\\
\end{eqnarray*}
respectively. It follows from a theorem of Lefshetz (see for example \cite{Lang}) that the $\BP^{15}$ embedding of the Jacobian is given by a basis of the space $\mathcal{L}(2(\Theta^{+} + \Theta^-))$ and the $\BP^3$ embedding of the Kummer surface is given by a basis of the space $\mathcal{L}(\Theta^{+} +\Theta^-)$. If $\Theta$ is the theta-divisor corresponding to any fixed $\mathfrak{k}$-rational Weierstrass point, where $\mathfrak{k}$ is an extension field of $k$, then $\mathcal{L}(\Theta^{+} +\Theta^-)$ is isomorphic to $\mathcal{L}(2\Theta)$ over $\mathfrak{k}$.

We do not give the $\BP^{15}$ embedding of the Jacobian here; it can be found, for example, in \cite{CasselsFlynn}. A $k$-rational point $P$ on $J$ can be represented as an unordered pair 
$\{ P_1,P_2\}$
where $P_1$ and $P_2$ are points on the curve $C$ that are either both defined over $k$ or are defined over a quadratic extension of $k$ and conjugate over $k$ such that $P_1+P_2-\infty^+-\infty^-$ or $P_1+P_2-2\infty$ is in $P$, viewed as a divisor class on $C$. If $P\ne 0$, then this representation is unique.
 
Following the notation from \cite{CasselsFlynn}, suppose $P_1=(x,y)$ and $P_2=(u,v)$ are affine points on the curve. Then a projective embedding of the Kummer surface is given by 
\begin{eqnarray*}
 \kappa_1&=&1\\
 \kappa_2&=&x+u\\
 \kappa_3&=&xu\\
 \kappa_4&=&\frac{F_0(x,u)-2yv}{(x-u)^2},\\
\end{eqnarray*}
where 
\begin{eqnarray*}
 F_0(x,u)&=&2f_0+f_1(x+u)+2f_2(xu)+f_3(x+u)xu+2f_4(xu)^2\\
         & &+f_5(x+u)xu+2f_6(xu)^3.\\
\end{eqnarray*}

The functions $\kappa_1,\kappa_2,\kappa_3,\kappa_4$ satisfy the quartic equation
\[
 K(\kappa_1,\kappa_2,\kappa_3,\kappa_4)=K_2(\kappa_1,\kappa_2,\kappa_3)\kappa^2_4+K_1(\kappa_1,\kappa_2,\kappa_3)\kappa_4+K_0(\kappa_1,\kappa_2,\kappa_3)=0,
\]
where
\begin{eqnarray*}
K_2(\kappa_1,\kappa_2,\kappa_3) &=& \kappa_2^2-4\kappa_1\kappa_3\\
K_1(\kappa_1,\kappa_2,\kappa_3) &=& -4\kappa_1^3f_0-2\kappa_1^2\kappa_2f_1-4\kappa_1^2\kappa_3f_2-2\kappa_1\kappa_2\kappa_3f_3-4\kappa_1\kappa_3^2f_4\\
																& &  -2\kappa_2\kappa_3^2f_5-4\kappa_3^3f_6\\
K_0(\kappa_1,\kappa_2,\kappa_3)  &=&-4\kappa_1^4f_0f_2+\kappa_1^4f_1^2-4\kappa_1^3\kappa_2f_0f_3-2\kappa_1^3\kappa_3f_1f_3-4\kappa_1^2\kappa_2^2
f_0f_4\\
&&+4\kappa_1^2\kappa_2\kappa_3f_0f_5-4\kappa_1^2\kappa_2\kappa_3f_1f_4-4\kappa_1^2\kappa_3^2f_0f_6+2\kappa_1^2\kappa_3^
2f_1f_5\\
&&-4\kappa_1^2\kappa_3^2f_2f_4+\kappa_1^2\kappa_3^2f_3^2-4\kappa_1\kappa_2^3f_0f_5+8\kappa_1\kappa_2^2\kappa_3f_0
f_6\\
&&-4\kappa_1\kappa_2^2\kappa_3f_1f_5+4\kappa_1\kappa_2\kappa_3^2f_1f_6-4\kappa_1\kappa_2\kappa_3^2f_2f_5\\
&&-2\kappa_1\kappa_3^3f_3
f_5-4\kappa_2^4f_0f_6-4\kappa_2^3\kappa_3f_1f_6-4\kappa_2^2\kappa_3^2f_2f_6\\
&&-4\kappa_2\kappa_3^3f_3f_6-4\kappa_3^4f_4f_6+\kappa_3^4f_5^2.\\
\end{eqnarray*}

We let $\kappa:=(\kappa_1,\kappa_2,\kappa_3,\kappa_4)$ be the map from the Jacobian into $\BP^3$. Clearly it identifies inverses and is hence $2:1$, except on points of order~2, where it is injective, so its image is an explicit realization of the Kummer surface $K$ in $\BP^3$ given by the defining equation $K(\kappa_1,\kappa_2,\kappa_3,\kappa_4)=0$. 

The question is how the group law is reflected on the Kummer surface. Firstly, since a point $Q\in J$ of order~2 is equal to its inverse and  $K$ precisely identifies inverses, addition of $\kappa(Q)$ is well-defined on the Kummer surface $K$. Furthermore, addition of $\kappa(Q)$ extends to a linear map on $\BP^3$, since it leaves $\left[\Theta^+ + \Theta^-\right]$ invariant, and thus can be given as multiplication by a matrix $W$. 

Secondly, there is a matrix $B=(B_{ij})_{i,j\in\{1,2,3,4\}}$ of biquadratic forms with the following property: Suppose $\overline{x}=(x_1,x_2,x_3,x_4)$ and $\overline{y}=(y_1,y_2,y_3,y_4)$ are quadruples such that $(x_1:x_2:x_3:x_4)=\kappa(P)$ and $(y_1:y_2:y_3:y_4)=\kappa(Q)$ for some $P,Q\in J$ (which we will abbreviate by saying that $\overline{x}$ and $\overline{y}$ are sets of Kummer coordinates for $P,Q$ respectively), then we can choose Kummer coordinates $\overline{w}$ and $\overline{z}$ for $P+Q$ and $P-Q$ respectively such that 
\[
 B_{ij}(\overline{x},\overline{y})=(w_i z_j + w_j z_i).
\]

Finally, multiplication by~2 is well-defined on the Kummer surface, because duplication commutes with negation. This duplication map is given by quartic polynomials $\delta_1,\delta_2,\delta_3,\delta_4$ (unique modulo the defining equation of the Kummer surface) such that for $P\in J$ we have  
\[
 \kappa(2P)=\delta(\kappa(P)),
\]
where $\delta=(\delta_1,\delta_2,\delta_3,\delta_4)$. Explicit expressions for $W,B$ and $\delta$ can be found in \cite{CasselsFlynn}.

The crucial point is that the existence of all of these objects can be asserted using quite general arguments (cf. \cite{Flynnb}) that do not involve the model of $C$, so we know that in the general situation we want to consider they can be found in principle. The question is how to actually compute them.

\section{Kummer coordinates for arbitrary characteristic} \strut\\[1mm] \label{GenKum}
We want to generalize Flynn's construction to the case of a genus~2 curve $C$ defined over a field $k$ of arbitrary characteristic, hence we need to consider affine defining equations of the form 
\[
 y^2+h(x)y=f(x),
\]
where
\[
f(x)=f_0+f_1x+f_2x^2+f_3x^3+f_4x^4+f_5x^5+f_6x^6
\] and 
\[
h(x)=h_0+h_1x+h_2x^2+h_3x^3
\]
are polynomials in $k[x]$. Let $C$ be such a curve and let $J$ denote its Jacobian.
The first obvious task is to find the map $\kappa:J\To \BP^3$ in the general case. 

As in \cite{Flynnb} we are required to find a basis for the~4-dimensional vector space 
\[
 \mathcal{L}(\Theta^+ + \Theta^-),
\]
since such a basis will give the desired map $\kappa:J\To \BP^3$. Suppose we have a generic point $P\in J$ represented by an unordered pair $\{P_1,P_2\}$, where $P_1=(x,y)$ and $P_2=(u,v)$. A basis may be found by looking for four linearly independent functions on $J$ which are symmetric in $P,Q$, have a pole of order at most~1 at infinity and may have a pole of any order at $0\in J$, but are regular elsewhere. Now, as in \cite{Flynnb}, 3 members of such a basis are easily found, namely the symmetric polynomials in $x$ and $u$ given by $\kappa_1=1,\kappa_2=x+u$ and $\kappa_3=xu$.

Looking for a suitable fourth coordinate, the following basis can be found:
\[
\kappa_1=1,\kappa_2=x+u,\kappa_3=xu,\kappa_4=\frac{F_0(x,u)-2yv-h(x)v-h(u)y}{(x-u)^2}.
\]
This obviously specializes to the basis given in section \ref{classical} in the case $h=0$ and it also specializes to the basis introduced in \cite{SD} when $\Char(k)=2$ and $h_3=0$. All of these are elements of $\mathcal{L}(\Theta^+ + \Theta^-)$, because they are even, symmetric, have no pole except at infinity, and grow at worst like $xu$ at infinity. We have a basis, because these~4 elements of the~4-dimensional vector space $\mathcal{L}(\Theta^+ + \Theta^-)$ are obviously linearly independent.

Similar to the classical case, these $\kappa_1,\kappa_2,\kappa_3,\kappa_4$ satisfy the quartic equation
\[
 K(\kappa_1,\kappa_2,\kappa_3,\kappa_4)=K_2(\kappa_1,\kappa_2,\kappa_3)\kappa^2_4+K_1(\kappa_1,\kappa_2,\kappa_3)\kappa_4+K_0(\kappa_1,\kappa_2,\kappa_3)=0,
\]
where
\begin{eqnarray*}
K_2(\kappa_1,\kappa_2,\kappa_3) &=& \kappa_2^2-4\kappa_1\kappa_3,\\
K_1(\kappa_1,\kappa_2,\kappa_3) &=& -4f_2\kappa_1^2\kappa_3-4f_6\kappa_3^3-4f_0\kappa_1^3
 -h_1h_3(\kappa_2^2\kappa_3-2\kappa_1\kappa_3^2)\\
 &&-h_2h_3\kappa_2\kappa_3^2-h_1h_2\kappa_1\kappa_2\kappa_3
 -h_1^2\kappa_1^2\kappa_3-2f_3\kappa_1\kappa_2\kappa_3-h_0^2\kappa_1^3\\
 &&-h_2^2\kappa_1\kappa_3^2 -2f_5\kappa_2\kappa_3^2-h_3^2\kappa_3^3-4f_4\kappa_1\kappa_3^2-2f_1\kappa_1^2\kappa_2\\
 &&-h_0h_1\kappa_1^2\kappa_2-h_0h_2(\kappa_1\kappa_2^2-2\kappa_1^2\kappa_3)-h_0h_3(\kappa_2^3-3\kappa_1\kappa_2\kappa_3),\\
K_0(\kappa_1,\kappa_2,\kappa_3)  &=& (-4f_0f_2 - f_0h_1^2 + f_1^2 + f_1h_0h_1 - f_2h_0^2)\kappa_1^4\\
 &&  + (-4f_0f_3 - 2f_0h_1h_2 + f_1h_0h_2 - f_3h_0^2)\kappa_1^3\kappa_2 \\
 &&  + (2f_0h_1h_3 - 2f_1f_3 - f_1h_0h_3 - f_1h_1h_2 + 2f_2h_0h_2 \\
 && \ \  - f_3h_0h_1)\kappa_1^3\kappa_3 \\
 && + (-4f_0f_4 - 2f_0h_1h_3 - f_0h_2^2 + f_1h_0h_3 - f_4h_0^2)\kappa_1^2\kappa_2^2 \\
 && + (4f_0f_5 + 2f_0h_2h_3 - 4f_1f_4 - f_1h_1h_3 - f_1h_2^2 + 2f_2h_0h_3 \\
 && \ \  + f_3h_0h_2 - 2f_4h_0h_1 + f_5h_0^2)\kappa_1^2\kappa_2\kappa_3 \\
 && + (-4f_0f_6 - f_0h_3^2 + 2f_1f_5 + f_1h_2h_3 - 4f_2f_4 - f_2h_2^2 + f_3^2\\
 && \ \  + f_3h_0h_3 + f_3h_1h_2 - f_4h_1^2 + f_5h_0h_1 - f_6h_0^2)\kappa_1^2\kappa_3^2\\
 &&  + (-4f_0f_5 - 2f_0h_2h_3 - f_5h_0^2)\kappa_1\kappa_2^3\\ 
 && + (8f_0f_6 + 2f_0h_3^2 - 4f_1f_5 - 2f_1h_2h_3 + f_3h_0h_3 - 2f_5h_0h_1\\
 && \ \  + 2f_6h_0^2)\kappa_1\kappa_2^2\kappa_3 \\
 && + (4f_1f_6 + f_1h_3^2 - 4f_2f_5 -  2f_2h_2h_3+ f_3h_1h_3\\
 && \ \  + 2f_4h_0h_3 - f_5h_0h_2 - f_5h_1^2 + 
      2f_6h_0h_1)\kappa_1\kappa_2\kappa_3^2\\
 && + (-2f_3f_5 - f_3h_2h_3 + 2f_4h_1h_3 - f_5h_0h_3 - 
    f_5h_1h_2\\
 && \ \  + 2f_6h_0h_2)\kappa_1\kappa_3^3  \\
 && + (-4f_0f_6 - f_0h_3^2 - f_6h_0^2)\kappa_2^4\\
 && + (-4f_1f_6 - f_1h_3^2 - 2f_6h_0h_1)\kappa_2^3\kappa_3 \\
 && + (-4f_2f_6 - f_2h_3^2 + f_5h_0h_3 - 2f_6h_0h_2 - f_6h_1^2)\kappa_2^2\kappa_3^2\\
 && + (-4f_3f_6 - f_3h_3^2 + f_5h_1h_3 - 2f_6h_1h_2)\kappa_2\kappa_3^3 \\
 && + (-4f_4f_6 - f_4h_3^2 + f_5^2 + f_5h_2h_3 - f_6h_2^2)\kappa_3^4.\\
\end{eqnarray*}
The zero locus of $K(\kappa_1,\kappa_2,\kappa_3,\kappa_4)$ gives an explicit realization of the Kummer surface associated with $J$. Notice that this equation becomes the classical one from \cite{Flynnb} in the case $C:y^2=f(x)$ and it reduces to the equation satisfied by the Kummer surface obtained by Duquesne in \cite{SD} if $\Char(k)=2$ and $h_3=0$.

Now our task is to compute the maps on the Kummer surface which make it so useful, namely the duplication map $\delta$, the matrix of biquadratic forms $B$ and the matrix $W$ that corresponds to translation by a point of order~2.

\section{Duplication}\strut\\[1mm] \label{Duplication}
We start by calculating the duplication map $\delta=(\delta_1,\delta_2,\delta_3,\delta_4)$ on the Kummer surface $K$. Here $\delta_1,\delta_2,\delta_3,\delta_4$ are quartic polynomials in $\kappa_1,\kappa_2,\kappa_3,\kappa_4$ and $\delta$ makes the following diagram commute:

\[
\xymatrix{J\ar[d]^{\kappa}\ar[r]^{[2]} & J \ar[d]^{\kappa} \\
	K \ar[r]^{\delta} & K,}
\]
where $[2]$ denotes the multiplication-by-2 map on the Jacobian.

In the classical case, Flynn uses the biquadratic forms described in the next section to obtain the duplication map. This is possible in the present situation (and gives the same result), but we can also use a different approach that does not depend on the biquadratic forms. We temporarily assume that $k$ is a field of characteristic not equal to~2, so that we can find a simpler model $C'$ for our curve $C$ given by
\[
 y^2=4f(x)+h(x)^2.
\]
Let $J'$ denote its Jacobian and let $K'$ denote its Kummer surface. Then clearly $J$ and $J'$ are isomorphic, as are $K$ and $K'$, so if we can explicitly determine the isomorphism
\[
 \tau:K\isom K'
\]
induced by the isomorphism $C\isom C'$, we can make use of the following commutative diagram, where $\delta'$ denotes the duplication map on $K'$:
\[
\xymatrix{K\ar[d]^{\tau}\ar[r]^{\delta} & K \ar[d]^{\tau} \\
	K' \ar[r]^{\delta'} & K'}
\]
It is easy to find the isomorphism $\tau$, in fact a short calculation shows that it is given by
\begin{eqnarray*}
	\tau: K&\To & K'\\
		 (\kappa_1,\kappa_2,\kappa_3,\kappa_4)&\mapsto& \left(\kappa_1,\kappa_2,\kappa_3,4\kappa_4-2(h_0h_2\kappa_1+h_0h_3\kappa_2+h_1h_2\kappa_3)\right).\\
\end{eqnarray*}
Thus we can find $\delta$ as
\[
 \delta:=\tau\circ\delta'\circ\tau^{-1}.
\]

Notice that this construction is only valid for characteristic $\ne 2$, so in order to remain valid in the remaining case, we want the polynomials $\delta_i$ to be defined and remain non-trivial modulo~2. Unfortunately this is not the case, but we can use the fact that the duplication map is only defined modulo the defining polynomial $K(\kappa_1,\kappa_2,\kappa_3,\kappa_4)$ and hence we can add multiples of this polynomial to the $\delta_i$. We do not change $\delta_1$ and $\delta_3$, but we add $-(32h_0h_3+32h_1h_2)\times K(\kappa_1,\kappa_2,\kappa_3,\kappa_4)$ to $\delta_2$ and $(48h_0h_1h_2h_3+48h_0^2h_3^2+32h_0h_3f_3)\times K(\kappa_1,\kappa_2,\kappa_3,\kappa_4)$ to $\delta_4$. After dividing all the $\delta_i$ by $64$ we obtain polynomials, also called $\delta_1,\delta_2,\delta_3,\delta_4$, that are defined and remain non-trivial modulo $2$.

\begin{prop}\label{propdelta}
The map $\delta$ constructed above represents duplication on the Kummer surface in any characteristic.
\end{prop}
\begin{proof}
We only need to show that the map $\delta=(\delta_1,\delta_2,\delta_3,\delta_4)$ represents duplication in characteristic~2. Since this is a geometric statement, we may as well assume that we have a field $k$ of characteristic~2 that is algebraically closed. Let $W(k)$ be its ring of Witt vectors with field of fractions $\mathfrak{k}$. Let 
\[
\widetilde{C}:y^2+\widetilde{h}(x)y=\widetilde{f}(x)
\]
be a genus~2 curve defined over $k$, with Jacobian $\widetilde{J}$ and Kummer surface $\widetilde{K}$ then $\widetilde{C}$ lifts to a genus~2 curve $C$ over $\mathfrak{k}$; $\widetilde{J}$ and $\widetilde{K}$ lift to the Jacobian $J$ and Kummer surface $K$ of $C$, repsectively.
 Then a Kummer surface $\widetilde{K}$ of a Jacobian $\widetilde{J}$ lifts to a Kummer surface $K$ of a Jacobian $J$ over $\mathfrak{k}$. Let $\delta$ denote the duplication map on $K$ that we have just found, reducing to the well-defined, non-trivial map $\widetilde{\delta}$ on $\widetilde{K}$. Let $\widetilde{P}\in \widetilde{J}$, lifting to $P\in J$. Then 
\[
 \delta(\kappa(P))=\kappa(2P)
\] 
and so if we normalize $\kappa(P)$ such that the entries lie in $W(k)$ with one of them having valuation zero, then either
\[
 \widetilde{\delta}(\widetilde{\kappa}(\widetilde{P}))= \widetilde{\kappa}(2 \widetilde{P})
\]
or $\widetilde{\delta_1}(\widetilde{\kappa}(\widetilde{P}))=\ldots=\widetilde{\delta_4}(\widetilde{\kappa}(\widetilde{P}))=0$. This can be seen by viewing $\delta$ and $\widetilde{\delta}$ as maps on the respective $\BP^3$'s. We need to show that the latter case cannot occur. 

For this we first reduce to a few simple cases. We use a transformation of the curve so that, depending on the number of roots of the homogenization $H(X,Z)$ of homogenous degree~3 of the polynomial $h(x)$, we are in one of the following three situations:
\begin{itemize}
 \item[(a)] $h=1$ 
 \item[(b)] $h=x$
 \item[(c)] $h=x^2+x$
\end{itemize}
Next, we can use another suitable transformation $y\mapsto y+u(x)$ where $u(x)$ is a polynomial of degree at most~3. It is not difficult to see that we can reduce to the case that 
\[
 f=f_1x+f_3x^3+f_5x^5
\]
where the condition that $C$ is nonsingular means in the respective cases:
\begin{itemize}
 \item[(a)] $f_5\ne0$
 \item[(b)] $f_1f_5\ne0$
 \item[(c)] $f_1f_5(f_1+f_3+f_5+f^2_1+f^2_3+f^2_5)\ne 0$
\end{itemize}
For each of these cases let $\overline{x}=(x_1,x_2,x_3,x_4)\in k^4$ be a quadruple that satisfies the defining equation $\widetilde{K}(\overline{x})=0$ of the Kummer surface associated to the Jacobian of $C$. We can use elementary methods, quite similar to those used to prove proposition 3.1 in \cite{StollH2}, to show the following.
\begin{lemma}\label{deltaallchar}
If $\widetilde{\delta_i}(\overline{x})=0$ for all $i\in\{1,2,3,4\}$, then we must already have $x_i=0$ for all $i\in\{1,2,3,4\}$.
\end{lemma}
 This means that the quadruple does not define a point on the Kummer surface and so the map $\widetilde{\delta}$ represents the duplication map on $K$. Since the proofs are not very enlightening but rather lengthy, they are not given here but may be found in section \ref{proofdelta}.
\end{proof}


In the special case $\Char(k)=2$ and $h_3=0$, the map $\delta$ specializes to the map given for duplication in \cite{SD} and in the case that the curve is given by an affine equation $y^2=f(x)$, it coincides with the duplication map given in \cite{Flynnb}.

\section{Biquadratic forms}\strut\\[1mm] \label{Biquads}
Let $P,Q\in J$ and let $\overline{x},\overline{y}$ be Kummer coordinates for $P$ and $Q$ respectively. The addition on the Jacobian does not descend to give a well-defined addition map on the Kummer surface. Indeed, given $\overline{x}$ and $\overline{y}$, we can find Kummer coordinates of $\kappa(P+Q)$ and $\kappa(P-Q)$, but in general we cannot tell them apart. Instead we can deduce from classical identities of theta-functions (see \cite{Hudson}) that projectively for any $i,j\in\{1,2,3,4\}$ $\kappa_i(P+Q)\kappa_j(P-Q)+\kappa_j(P+Q)\kappa_i(P-Q)$ is biquadratic in the $(x_1,x_2,x_3,x_4),(y_1,y_2,y_3,y_4)$ and therefore there is a matrix $B:=(B_{ij})_{i,j\in\{1,2,3,4\}}$ of biquadratic forms in $x,y$ having the property that there are Kummer coordinates $\overline{w}$ and $\overline{z}$ for $P+Q$ and $P-Q$ respectively such that 
\[
B_{ij}(\overline{x},\overline{y})=w_i z_j + w_j z_i.
\]
For its computation, we will again use the fact that the Kummer surface $K$ is isomorphic to $K'$ defined in the last section. The isomorphism $\tau: K\To K'$ was also given there. 

Let $B'$ denote the corresponding matrix of biquadratic forms on $K'$ and let $\overline{x}'=\tau(\overline{x}),\overline{y}'=\tau(\overline{y}),\overline{z}'=\tau(\overline{z}), \overline{w}'=\tau(\overline{w})$, so that we have
\begin{equation}\label{Bprime}
 B'_{i,j}(\overline{x}',\overline{y}')=w'_i z'_j + w'_j z'_i.
\end{equation}
Notice that for $i\in\{1,2,3\}$, we have $x'_i=x_i,y'_i=y_i,z'_i=z_i,w'_i=w_i$. We use this fact, our explicit expression of the isomorphism $\tau$ and \eqref{Bprime} to find the matrix $B$ in terms of the entries of $B'$. We write $b'_{i,j}$ for $B'_{i,j}(\overline{x}',\overline{y}')$.

For $i,j\in\{1,2,3\}$ we have
\[
 B_{i,j}(\overline{x},\overline{y})=w_i z_j + w_j z_i=w'_i z'_j + w'_j z'_i=b'_{i,j}.
\]
To find an entry of the fourth column (or row) of $B$ not equal to $b'_{4,4}$ we have to do some algebra. We get, for example,
\[
 B_{1,4}(\overline{x},\overline{y})=\frac{1}{4}b'_{1,4}+\frac{1}{2}\left(2h_0h_2b'_{1,1}+h_0h_3b'_{1,2}+h_1h_3b'_{1,3}\right)
\]
and analogous formulas for $B_{2,4}(\overline{x},\overline{y})$ and $B_{3,4}(\overline{x},\overline{y})$. Finally we compute
\[ B_{4,4}(\overline{x},\overline{y})=\frac{1}{4}\left(h_0h_2b'_{1,4}+h_0h_3b'_{2,4}+h_1h_3b'_{3,4}+h^2_0h^2_2b'_{1,1}+h^2_0h^2_3b'_{2,2}+h^2_1h^2_3b'_{3,3}\right)\]\[ + \frac{1}{8}\left(h^2_0h_2h_3b'_{1,2}+h_0h_1h_2h_3b'_{1,3}+h_0h_1h^2_3b'_{2,3}\right)+\frac{1}{16}b'_{4,4}. 
\]
Dividing all~16 entries of the matrix thus computed by~16, we obtain a matrix whose entries are all defined and remain non-trivial modulo~2. 

\begin{prop}
We have $B_{i,j}(\overline{x},\overline{y})=w_i z_j + w_j z_i$ in any characteristic.
\end{prop}
\begin{proof}
As in the last section, we are required to verify that this matrix actually contains the biquadratic forms we were looking for in characteristic~2. Keeping the notation from section \ref{Duplication}, we let $\widetilde{B_{i,j}}$ denote the reduction of the biquadratic form $B_{i,j}$ on a Kummer surface $K$ over the fraction field of the ring of Witt vectors reducing to $\widetilde{K}$. Viewing the $B_{i,j}$ and the $\widetilde{B_{i,j}}$ as maps on $\BP_\mathfrak{k}^3\times\BP_\mathfrak{k}^3$ and $\BP_k^3\times\BP_k^3$ respectively, we see that for a given point $((x_1:x_2:x_3:x_4),(y_1:y_2:y_3:y_4))\in \widetilde{K}^3\times \widetilde{K}^3$ either all $\widetilde{B_{i,j}}(x_1,x_2,x_3,x_4;y_1,y_2,y_3,y_4)$ vanish or they give the correct biquadratic forms. 

The proof of the proposition is finished by the following lemma:
\begin{lemma}\label{Ballchar}
If $\overline{x}=(x_1,x_2,x_3,x_4)\in k^4$ and $\overline{y}=(y_1,y_2,y_3,y_4)\in k^4$ satisfy $\widetilde{K}(\overline{x})=\widetilde{K}(\overline{y})=0$ and if $\widetilde{B_{i,j}}(\overline{x},\overline{y})$ all vanish, then $x_i=0$ for all $i$ or $y_i=0$ for all $i$.
\end{lemma}
By the discussion in section \ref{Duplication} we can reduce to the cases (a), (b) and (c) introduced there. The proofs for these cases can be found in section \ref{proofB}. Note that the methods are again similar to those employed in the proof of proposition 2.1 of \cite{StollH2}; they consist of straightforward, but quite lengthy, algebraic manipulations.
\end{proof}

In fact, in the case of characteristic~2 the matrix $B$ reduces to the corresponding matrix found in \cite{SD} when $h_3=0$ and in the case $h(x)=0$ it reduces to the matrix worked out by Flynn in \cite{Flynnb}.

\section{Translation by a point of order~2}\strut\\[1mm] \label{Wg}
Let $Q\in J$ be a point of order~2, so $P + Q = P - Q$ for all $P\in J$ and translation by $\kappa(Q)$ is defined on the Kummer surface. In fact, it is a linear map on $\BP^3$, so it can be given as a matrix in terms of the coefficients of the curve. This matrix was found in the special case $C:y^2=f(x)$ by Flynn in \cite{Flynnb} and is given in terms of the coefficients of polynomials $s$ and $t$, where $f(x)=s(x)t(x)$, $\deg(s)=2, \deg(t)=4$ and the roots of $s$ are the $x$-coordinates of the points $Q_1,Q_2$ on the curve $C$ such that $Q$ can be represented by the unordered pair $\{Q_1,Q_2\}$. Furthermore, the map is an involution and hence the square of the matrix representing it is a scalar multiple of the identity matrix.

As before, we proceed by making use of the isomorphism $\tau:K\To K'$ in the case $\Char(k)\ne 2$. Let $W'$ denote the matrix corresponding to translation by $\tau(\kappa(Q))$ on $K'$. We want to find the matrix $W$ that makes the following diagram commute
\[
\xymatrix{K\ar[d]^{\tau}\ar[r]^{W} & K \ar[d]^{\tau} \\
	K' \ar[r]^{W'} & K',}
\]
where the horizontal maps are multiplication by the respective matrix. This means that we will express the resulting matrix in terms of the coefficients of polynomials $s,t$ such that $4f(x)+h(x)^2=s(x)t(x)$. First we compute
\[
 W:=T^{-1}W'T,
\]
where $T$ is the matrix corresponding to $\tau$. Then $W$ has the desired properties for $\Char(k)\ne 2$. In order to generalize it to arbitrary characteristic, one could try to manipulate the entries directly, or one could first express them in terms of the Kummer coordinates of $Q$, as opposed to the coefficients of $s$ and $t$. Unfortunately, neither of these approaches has proved successful, see the discussion below. Therefore, we subsequently use a different method to compute the matrix corresponding to translation by a point of order~2 when $\Char(k)=2$. Our method is analogous to the one used by Flynn in the case that $\Char(k)\ne 2$ and $h=0$. In addition, it is identical with the method used independently by Duquesne in the case that $\Char(k)=2$ and $h$ has degree~2. However, the matrix computed there only works when $Q$ does not involve a point at infinity.

Suppose that $C$ is a curve of genus~2 given by an affine equation $C:y^2+h(x)y=f(x)$ and defined over a field $k$ of characteristic equal to~2. Let $Q$ be a $k$-rational point of order~2 on its Jacobian $J$. In order to find the matrix $W$ corresponding to translation by $Q$, we directly compute the image of $P+Q$ on the Kummer surface using the geometric group law on the Jacobian, where $P\in J(k)$ is generic, and then simplify to make it linear in the Kummer coordinates of $P$. 
The point $Q$ can be represented as $\{Q_1,Q_2\}$ with points $Q_i\in C$. First we assume that $Q_1$ and $Q_2$ are affine points, so $Q_i=(x_i,y_i)$ satisfying 
\[
 h(x)=(x-x_1)(x-x_2)t(x),
\]
where $t(x)=t_0+t_1x$. 

We will keep the discussion of this case brief (see \cite{SD} or \cite{Flynnb} for a more detailed discussion). We start by finding the first~3 rows of the matrix $W$ such that $W\kappa(P)=\kappa(P+Q)$; the last row is computed using the fact that $W^2$ must be a scalar multiple of the identity matrix. 
After a little simplification the matrix can be expressed in terms of the Kummer coordinates $k_1,k_2,k_3,k_4$ of $Q$ and the coefficients of the polynomials $f,t$ and $b$, where $y=b(x)=b_1+b_0x$ is the line joining the points $Q_1$ and $Q_2$, so
\[
   b_0=\frac{y_1-y_2}{x_1-x_2},\;\; b_1=\frac{x_2y_1-x_1y_2}{x_1-x_2}.
\]
Recall that a point on the Jacobian can be given in Mumford representation as $(a(x),b(x))$, where $a(x)=(x-x_1)(x-x_2)=x^2-\frac{k_2}{k_1}x+\frac{k_3}{k_1}$. 
 
To complete the picture, we have to find the matrix $W$ in the case that $Q_1=(x_1,y_1)$ is affine and $Q_2$ is at infinity. Then $b(x)$ is a cubic polynomial. Its leading coefficient $r_6$ plays the role of the $y$-coordinate of $Q_2$ and we distinguish between the cases $Q_2=\infty^+$ and $Q_2=\infty^-$ according to the value of $r_6$. By going through the same steps as before, we find $W$ in terms of $r_6,y_1$, the coefficients of $f$ and $t$ and the Kummer coordinates of $Q$.

In order to unify the two matrices, the following notation is convenient:
We set $k'_i:=k_i/k_2$ in both cases. If $Q_2$ is affine we set
\begin{eqnarray*}
b'_0&:=&\frac{y_1-y_2}{(x_1-x_2)^2}=\frac{b_0}{x_1-x_2},\\
b'_1&:=&\frac{y_1x_2-y_2x_1}{(x_1-x_2)^2}=\frac{b_1}{x_1-x_2},\\
b'_2&:=&\frac{y_1x_2^2-y_2x_1^2}{(x_1-x_2)^2}=b'_1\frac{k'_2}{k'_1}+b'_0\frac{k'_3}{k'_1},\\
b'_3&:=&\frac{y_1x_2^3-y_2x_1^3}{(x_1-x_2)^2}=b'_2\frac{k'_2}{k'_1}+b'_1\frac{k'_3}{k'_1}=b'_1\left(\frac{k'_2}{k'_1}\right)^2+b'_1\frac{k'_3}{k'_1}+b'_0\frac{k'_2k'_3}{k'^2_1},\\
c&:=&\frac{y_1y_2}{x_1-x_2}=b'_0b'_1\left(\frac{k'_2}{k'_1}\right)^3+\frac{f(x_1)x_2+f(x_2)x_1}{x_1-x_2}.
\end{eqnarray*} 

Now suppose that $Q_2$ is at infinity. In this situation we set
\begin{eqnarray*}
 b'_i&:=&r_6k'^i_3 \text{  for }i=0,1,2,\\
 b'_3&:=&r_6k'^3_3+y_1, \\
 c&:=&y_1r_6.
\end{eqnarray*}
Here $y_1$ satisfies $y^2_1=f(x_1)$, hence it can be computed using the coefficients of $f$ and the $k'_i$, or by $y_1=b(x_1)$.

Then the unified matrix is given by
\[
W=\left( \begin{array}{cccc}
t_1b'_2+k'_4 & t_1b'_1+f_5k'_3 & t_1b'_0+f_5k'_2 & k'_1 \\

t_0b'_2+t_1b'_3+f_3k'_3 & t_0b'_1+t_1b'_2+k'_4 & t_0b'_0+t_1b'_1+f_3k'_1 & k'_2 \\

t_0b'_3+f_1k'_2 & t_0b'_2+f_1k'_1 & t_0b'_1+k'_4 & k'_3\\

W_{4,1}& W_{4,2} & W_{4,3} & k'_4
\end{array} \right),
\]
where
\begin{eqnarray*}
W_{4,1}&=&t_0f_1b'_0+t_0f_3b'_2+t_0^2c+t_1f_1b'_1+f_3f_1k'_1, \\
W_{4,2}&=&t_0f_5b'_3+t_0t_1c+t_1f_1b'_0+f_1f_5k'_2,\\
W_{4,3}&=& t_0f_5b'_2+t_1f_3b'_1+t_1f_5b'_3+t_1^2c+f_3f_5k'_3.
\end{eqnarray*}

It seems curious that our results in this section apparently cannot be combined to form a matrix that works in arbitrary characteristic. One possible reason for this is the fact that if $\Char(k)=2$, then an affine point $(x,y)$ invariant under the hyperelliptic involution satisfies $h(x)=0$ and if $\Char(k)\ne 2$, then such a point satisfies $y=0$. In general, we can only assume that $2y+h(x)=0$ and this is not a sufficient simplification to make the method used above work. Moreover, if $\Char(k)=2$, then, depending on the number of distinct roots of $h$, we have $\#J[2]\in \{1,2,4\}$, whereas otherwise $\#J[2]=16$. It would be interesting to find out whether there is a matrix $W$ representing translation by a point of order~2 in arbitrary characteristic, either by finding such a matrix or by proving that it cannot exist.

\section{Proof of lemma \ref{deltaallchar}}\label{proofdelta}

In this section we prove lemma \ref{deltaallchar} using case distinctions and elementary algebraic manipulations. It would be interesting to find a more conceptual proof.In all cases $f$ is of the form $f=f_1x+f_3x^3+f_5x^5$.  
\subsection{Case (a): $h=1$, $f_5\ne0$}
$\widetilde{\delta_i}(\overline{x})=\widetilde{K}(\overline{x})=0$ implies that
\[
 0=\widetilde{\delta_2}(\overline{x})+f_3\widetilde{K}(\overline{x})=f_5x^4_1,
\]
so $x_1=0$. We find $0=\widetilde{\delta_1}(\overline{x})=f_5x^4_2$ and hence $x_2=0$. Then we also obtain $x_3=0$ from $0=\widetilde{K}(\overline{x})=f_5x^4_3$ and thus $0=\widetilde{\delta_4}(\overline{x})=x^4_4$ means that indeed $x_i=0$ for all $i\in\{1,2,3,4\}$.
\subsection{Case (b): $h=x, f_1f_5\ne0$}
Similar to case (a) we have
\[
 0=\widetilde{\delta_2}(\overline{x})+\widetilde{\delta_3}(\overline{x})\widetilde{K}(\overline{x})=f_5x^2_1x^2_3,
\]
hence $x_1=0$ or $x_3=0$.

If $x_1=0$, then $0=\delta_3(\overline{x})=f_5x^4_3$, thus $x_3=0$. The Kummer surface equation then reads $0=K(\overline{x})=x^2_2x^2_4$, whence $x_2=0$ or $x_4=0$. However, if $x_2=0$, then $0=\delta_4(\overline{x})=x_4^4$ and if $x_4=0$, then $0=\delta_4(\overline{x})=f^2_1f^2_5x^4_2$. Therefore we can deduce that $x_i=0$ follows for all $i\in\{1,2,3,4\}$ in both subcases.

In the other case $x_3=0$ implies $0=\delta_3(\overline{x})=f^2_1x^4_1$, so $x_1=0$ and we are again in the situation already considered above.

\subsection{Case (c): $h=x^2+x,f_1f_5(f_1+f_3+f_5+f_1^2+f_3^2+f_5^2)\ne0$}
For this case, which is slightly more complicated than the two previous cases, we employ a case distinction on $x_1$. First we assume that $x_1=0$ and show that necessarily the other $x_i$ must be equal to zero as well. Then we suppose that $x_1\ne0$ and derive a contradiction. Furthermore, we abbreviate $\beta=f_1+f_3+f_5+f_1^2+f_3^2+f_5^2$.

So let $x_1=0$. Then
\[
 0=\delta_3(\overline{x})=x_3^2(x_4+f_5x_3)
\]
which means that we must have $x_3=0 \text{ or }x_4=f_5x_3$.

If $x_3=0$, then $0=K(\overline{x})=x_2^2x_4^2$ implies that $x_2=0$ or $x_4=0$. But from
\[ 
0=\delta_4(\overline{x})=x_4^4+f_1^2f_5^2x_2^4
\]
the result follows.

If we have $x_4=f_5x_3\ne 0$ instead, then 
\[
 0=K(\overline{x})=f_5^2x_3^2(x_2+x_3),
\]
so that we get $x_2=x_3\ne 0$ and hence
\[
 0=\delta_4(\overline{x})=f_5^2\beta x_3^4,
\]
a contradiction. This means that $x_1=0=\delta_i(\overline{x})$ is only possible if $x_i=0$ for all $i$.

Now we consider the case $x_1\ne0$, so we may assume that $x_1=1$. Here
\[
 0=\delta_2(\overline{x})+f_3K(\overline{x})=(1+x_2+x_3)x_3(x_4+f_1+f_5x_3).
\]
If $x_3=0$, we find that
\[
 0=\delta_1(\overline{x})=f_1^2+x_4^2\Rightarrow x_4=f_1\Rightarrow
K(\overline{x})=f_1^2(1+x_2)^2\Rightarrow x_2=1
\]
and so 
\begin{equation}\label{dx30}
 0=\delta_4(\overline{x})=f_1^2\beta,
\end{equation}
a contradiction.

Next we suppose that $x_4=f_1+f_5x_3$, leading to $0=\delta_1(\overline{x})=f_5^2x_3^2(1+x_2+x_3)^2$, so that either $x_3=0$ which leads to a contradiction by \eqref{dx30} or $1+x_2+x_3=0$ must hold. However, in that case we deduce that $0=K(\overline{x})=x_3^2\beta$, so we get a contradiction anyway.

Finally, we assume that $1+x_2+x_3=0$ and see that
\[
 0=\delta_1(\overline{x})=(x_4+f_1+f_5x_3)^2=0\Rightarrow x_4+f_1+f_5x_3=0,
\]
proving the lemma.

\section{Proof of lemma \ref{Ballchar}}\label{proofB}
The following section consists of a proof of lemma \ref{Ballchar}. In all cases our method is to first assume that $x_1=0$ and then show that either $x_i=0$ for all $i$ or $y_i=0$ for all $i$ follows. To finish the claim, we assume that $x_1\ne0$, so without loss of generality $x_1=1$, and then show that all $y_i$ must be zero. We abbreviate $\overline{x}=(x_1,x_2,x_3,x_4)$ and $\overline{y}=(y_1,y_2,y_3,y_4)$. There are a lot of nested case distinctions, so in order to follow the proof, the main difficulty is to remember at each step which assumptions were made. As in the case of lemma \ref{deltaallchar} a conceptual proof would be of interest.

\subsection{Case (a): $h=1$}
 First we assume that $x_1=0$. Then 
 \[
  0=B_{12}(\overline{x},\overline{y})=f_5x_2^2y_1^2\Rightarrow x_2=0 \text{ or } y_1=0
 \]
 If $x_2=0$, but $y_1\ne0$, then 
 \begin{eqnarray*}
 && 0=B_{14}(\overline{x},\overline{y})=f_5x_3^2y_1^2\Rightarrow x_3=0\\ &\Rightarrow & 0=B_{11}(\overline{x},\overline{y})=x_4^2y_1^2\Rightarrow x_4=0.
 \end{eqnarray*}
  If we have $y_1=0\ne x_2$ instead, then
 \begin{eqnarray*}
  &&0=B_{14}(\overline{x},\overline{y})=f_5^2x_2^2y_2^2  \Rightarrow  y_2=0\\ &\Rightarrow& 0=B_{22}(\overline{x},\overline{y})=x_2^2y_4^2\Rightarrow y_4=0\text{ and } \\
  && 0=B_{11}(\overline{x},\overline{y})=f_5^2x_2^2y_3^2\Rightarrow y_3=0.
 \end{eqnarray*}

The third case we have to look at is the case $x_2=y_1=0$. In this situation we get $0=B_{22}(\overline{x},\overline{y})=x_4^2y_2^2$, so $x_4=0$ or $y_2=0$. We also see that $0=K(\overline{x})=f_5^2x_3^4$ and hence $x_3=0$. So we may assume that $y_2=0\ne0$ which implies
\[
 0=B_{33}(\overline{x},\overline{y})=(y_3x_4)^2\text{ and }0=B_{44}(\overline{x},\overline{y})=x_4y_4,
\]
so $y_3=y_4=0$.

Now that we have finished proving that $x_1=0$ implies the lemma in case (a), the remaining step is to deduce that all $y_i$ must be qual to zero using the assumption that $x_1=1$. It follows quickly from the following observation
\[
 0=B_{12}(\overline{x},\overline{y})=f_5(y_2+x_2y_1)^2\Rightarrow y_2=x_2y_1,
\]
since then 
\begin{eqnarray*}
  &&0=B_{23}(\overline{x},\overline{y})=f_5(y_3+x_3y_1)^2  \Rightarrow  y_3=x_3y_1\\ &\Rightarrow& 0=B_{11}(\overline{x},\overline{y})=(y_4+x_4y_1)^2\Rightarrow y_4=x_4y_1 \\
  &\Rightarrow& 0=B_{24}(\overline{x},\overline{y})=f_5^2y_1^2\Rightarrow y_1=0\Rightarrow y_2=y_3=y_4=0.
 \end{eqnarray*}
\subsection{Case (b): $h=x$}
Suppose that $x_1=0$ and observe that $0=B_{12}(\overline{x},\overline{y})=f_5x_3^2y_1^2$, implying either $x_3=0$ or $y_1=0$.

If $x_3=0$, then we get
\begin{equation}\label{x30K}
0=K(\overline{x})=x_2^2x_4^2, \text{ so } x_2=0 \text{ or }x_4=0.
\end{equation}
If $x_2=0$, then 
\[
 0=B_{11}(\overline{x},\overline{y})=x_4^2y_1^2=B_{22}(\overline{x},\overline{y})=x_4^2y_2^2=B_{33}(\overline{x},\overline{y})=x_4^2y_3^2
\]
from which $x_i=0$ for all $i$ or $y_i=0$ for all $i$ follows. 

If $x_4=0$, then 
\begin{eqnarray*}
0&=&B_{11}(\overline{x},\overline{y})=f_5^2x_2^2y_3^2=B_{22}(\overline{x},\overline{y})=x_2^2y_4^2=B_{33}(\overline{x},\overline{y}=f_1^2x_2^2y_1^2\\
& =& B_{44}(\overline{x},\overline{y})=f_1^2f_5^2x_2^2y_2^2
\end{eqnarray*}
and so we find again that $x_i=0$ for all $i$ or $y_i=0$ for all $i$. 

Now we go back to \eqref{x30K} and suppose that $y_1=0\ne x_3$. However, this has the following consequence:
\[
 0=B_{34}(\overline{x},\overline{y})=f_5^2x_3^2y_3^2\Rightarrow y_3=0
\]
So we get 
$0=B_{33}(\overline{x},\overline{y})=x_3^2y_4^2\Rightarrow y_4=0$ and $0=B_{11}(\overline{x},\overline{y})=f_5^2x_3^2y_2^2$, therefore $y_2=0$.

We now consider the case $x_1=1$. Then 
\begin{equation}\label{y1f1}
0=B_{12}(\overline{x},\overline{y})=f_5(y_3+x_3y_1)^2\Rightarrow y_3=x_3y_1\Rightarrow 0=B_{34}(\overline{x},\overline{y})=y_1^2(f_1+f_5x_3^2)^2
\end{equation}
and hence either $y_1=0$ or we can express $f_1$ as $f_1=f_5x_3^2$.

The first case is $y_1=0$, which implies $y_3=0$ and $0=B_{23}(\overline{x},\overline{y})=x_3y_2y_4$.

If $x_3=0$, we get 
\[0= B_{11}(\overline{x},\overline{y})=y_4^2 \text{ and }0=B_{33}(\overline{x},\overline{y})=f_1^2y_2^2,\]
thus $y_2=y_4=0.$

If we have $y_2=0$ in \eqref{y1f1}, then again $0=B_{11}(\overline{x},\overline{y})=y_4^2\Rightarrow y_4=0$. Finally, if we have $y_4=0$ in \eqref{y1f1}, then $y_2=0$ since $0=B_{33}(\overline{x},\overline{y})=f_1^2y_2^2.$

In order to prove the lemma in case (b), it remains to prove it in the case $x_1=1, y_3=x_3y_1, f_1=f_5x_3^2$. This implies that $x_3\ne 0$ and hence $y_3=0$. We also obtain
\begin{equation}\label{y4long}
0=B_{33}(\overline{x},\overline{y})=x_3^2(y_4+x_4y_1+f_5x_3(y_2+y_1x_2))^2=0,
\end{equation}
whence $y_4=x_4y_1+f_5x_3(y_2+y_1x_2)$.

Using this relation we find that
\[
 0=B_{23}(\overline{x},\overline{y})=f_5x_3^2(y_2+x_2y_1)^2\Rightarrow y_2=x_2y_1
\]
and hence 
\[
0=B_{24}(\overline{x},\overline{y})=f_5x_3^2y_1^2\Rightarrow y_1=0\Rightarrow y_2=0.
\]
We also have $y_3=0$ from \eqref{y1f1} and $y_4=0$ because of \eqref{y4long}, which proves part (b) of the lemma.

\subsection{Case (c): $h=x^2+x, f_1f_5(f_1+f_3+f_5+f_1^2+f_3^2+f_5^2)\ne0$ }
Let $\beta=f_1+f_3+f_5+f_1^2+f_3^2+f_5^2$. This is the trickiest case of the lemma, although it is, like the other cases, completely elementary. Having said that, we again start off by assuming that $x_1=0$, which yields the following Kummer surface equation
\[
 0=K(\overline{x})=(f_5x_3^2+x_2x_4)^2\Rightarrow f_5x_3^2=x_2x_4
\]
which in turn implies
\begin{equation}\label{y1x3y3}
 0=B_{13}(\overline{x},\overline{y})=x_3y_1y_3(x_4+f_5x_2)\Rightarrow x_3=0 \text{ or } y_1=0\text{ or } y_3=0\text{ or } x_4=f_5x_2.
\end{equation}
We first assume that $y_1=0$ and get
\[
 0=B_{34}(\overline{x},\overline{y})=x_4y_3(y_4+f_5y_3)(x_2+x_3)
\]
and therefore
\begin{equation}\label{y3x4}
  x_4=0 \text{ or } y_3=0\text{ or } y_4=f_5y_3\text{ or } x_2=x_3.
\end{equation}

We will actually go through all the cases in \eqref{y3x4}. This is a rather tedious task, but we will be able to reuse several of the results in the other cases appearing in \eqref{y1x3y3}.

Suppose $y_3=0$. Then 
\begin{equation}\label{y3x3y4}
 0=B_{33}(\overline{x},\overline{y})=x_3^2y_4^2\Rightarrow x_3=0 \text{ or } y_4=0.
\end{equation}
If $x_3=0$, then 
\[
 0=K(\overline{x})=x_2^2x_4^2\Rightarrow x_2=0\text{ or }x_4=0.
\]
Now we get
\[
x_2=0\Rightarrow 0=B_{22}(\overline{x},\overline{y})=x_4^2y_2^2,0=B_{44}(\overline{x},\overline{y})= x_4^2y_4^2
\]
and
\[
x_4=0\Rightarrow 0=B_{22}(\overline{x},\overline{y})=x_2^2y_4^2,0=B_{44}(\overline{x},\overline{y})=f_1^2f_5^2 x_2^2y_2^2,
\]
so we see that in both cases either $x_i=0$ for all $i$ or $y_i=0$ for all $i$.

If $y_4=0$ in \eqref{y3x3y4}, then 
\[
0=B_{22}(\overline{x},\overline{y})=x_4^2y_2^2=B_{44}(\overline{x},\overline{y})=f_1^2f_5^2x_2^2y_2^2
\]
so either $y_2=0$ or $x_2=x_4=0$, in which case we have $0=K(\overline{x})=f_5^2x_3^4$. This finishes the case $y_3=0$ in \eqref{y3x4}.

If $x_4=0$ in \eqref{y3x4}, then 
\[
 0=B_{22}(\overline{x},\overline{y})=x_2^2y_4^2=B_{33}(\overline{x},\overline{y})=x_3^2y_4^2\Rightarrow y_4=0 \text{ or }x_2=x_3=0.
\]
Now $x_2=x_3=0$ means we are already done. If instead we have $y_4=0$, we get that
\[
0=B_{11}(\overline{x},\overline{y})=f_5^2x_2^2y_3^2\Rightarrow x_2=0\text{ or } y_3=0.
\]
We have already dealt with the case $y_3=0$, so we can assume $x_2=0$. But then we have $0=K(\overline{x})=f_5^2x_3^4$ again.

The next case in \eqref{y3x4} that we consider is the case $y_4=f_5y_3$ which implies $K(\overline{y})=f_5^2y_3^2(y_2+y_3)^2$. Since we know that $y_3=0$ implies our claim, we can assume that $y_2=y_3\ne0$. Then 
\begin{eqnarray*}
0&=&B_{33}(\overline{x},\overline{y})=y_3^2(x_4+f_5x_3)^2\\&&\Rightarrow x_4=f_5x_3\text{ and } 0=B_{22}(\overline{x},\overline{y})=f_5^2y_3^2(x_2+x_3)^2\\
&&\Rightarrow x_2=x_3
\end{eqnarray*}
and hence
\[
 0=B_{44}(\overline{x},\overline{y})=f_5^2\beta x_3^2y_3^2,
\]
so that finally $x_2=x_3=0$.

In order to finish off \eqref{y3x4} we assume that $x_2=x_3$, thus
\[
K(\overline{x})=x_3^2(f_5x_3+x_4)^2\Rightarrow x_3=0\text{ or }f_5x_3=x_4\ne0
\]
Assuming that $x_3=0$, we deduce from $x_2=x_3=0$ and
\[
0=B_{22}(\overline{x},\overline{y})=x_4^2y_2^2=B_{33}(\overline{x},\overline{y})=x_4^2y_3^2=B_{44}(\overline{x},\overline{y})=x_4^2y_4^2
\]
that either $x_i=0$ for all $i$ or $y_i=0$ for all $i$.

So we consider the case $x_4=f_5x_3\ne 0$ and see that
\[
 0=B_{33}(\overline{x},\overline{y})=x_3^2(y_4+f_5y_3)^2\Rightarrow y_4=f_5y_3
\]
Hence we have 
\[
 0=B_{22}(\overline{x},\overline{y})=f_5^2x_3^2(y_2+y_3)^2\Rightarrow y_2=y_3.
\]
But if $y_2=y_3$, then $0=B_{44}(\overline{x},\overline{y})=f_5^2\beta x_3^2y_3^2$ and so $y_3=0$, a case we have finished already. Therefore we have proved the assertion of the lemma for the case $x_1=y_1=0$.

Now we go back to \eqref{y1x3y3} and assume that $x_3=0$. The Kummer surface equation then tells us that either $x_2=0$ or $x_4=0$. But
\begin{eqnarray*}
 x_2=0\Rightarrow 0&=&B_{11}(\overline{x},\overline{y})=x_4^2y_1^2=B_{22}(\overline{x},\overline{y})=x_4^2y_2^2\\&=&B_{33}(\overline{x},\overline{y})=x_4^2y_3^2=B_{44}(\overline{x},\overline{y})=x_4^2y_4^2
\end{eqnarray*}
and
\begin{eqnarray*}
 x_4=0\Rightarrow 0&=&B_{11}(\overline{x},\overline{y})=f_5^2x_2^1y_3^2=B_{22}(\overline{x},\overline{y})=x_2^2y_4^2\\&=&B_{33}(\overline{x},\overline{y})=f_1^2x_2^2y_1^2=B_{44}(\overline{x},\overline{y})=f_1^2f_5^2x_2^2y_2^2,
\end{eqnarray*}
thus we get that in both cases either $x_i=0$ for all $i$ or $y_i=0$ for all $i$.

The next possible case from \eqref{y1x3y3} is $y_3=0$. Because of what we have shown already, we can assume that $y_1x_3\ne0$. We find that
\[
0=B_{23}(\overline{x},\overline{y})=y_1x_3(x_2+x_3)(f_1y_1+y_4),
\]
so that either $x_2=x_3\ne0$ or $y_4=f_1y_1\ne0$. 
In the former case we have $0=B_{33}(\overline{x},\overline{y})=x_3^2(y_4+f_1y_1)^2$, so we are in the latter case anyway.

Accordingly we suppose that $y_4=f_1y_1\ne0$ which means 
\[
 K(\overline{y})=f_1^2y_1^2(y_1+y_2)^2.
\]
Thus $y_2=y_1\ne0$ and from $0=B_{44}(\overline{x},\overline{y})=f_1^2y_1^2(x_4+f_5x_2)^2$ we get $x_4=f_5x_2$ which ultimately leads to
\[
 0=B_{22}(\overline{x},\overline{y})=\beta x_3^2y_1^2,
\]
a contradiction. This finishes case (c) of the lemma in the case $x_1=0$.

Now we assume that $x_1=1$. It turns out that it is a good idea to further distinguish between the cases $y_3=0$ and $y_3\ne0$.

We start with the case $y_3=0$ which leads to 
\begin{equation}\label{y4b11}
0=B_{11}(\overline{x},\overline{y})=(y_4+x_4y_1+f_5x_3y_2)^2\Rightarrow y_4+x_4y_1+f_5x_3y_2=0
\end{equation}
and thus 
\begin{equation}\label{y1y2b12}
0=B_{12}(\overline{x},\overline{y})=f_5x_3(y_1+y_2)^2=0
\end{equation}
so that either $x_3=0$ or $y_1=y_2$.

The assumption $x_3=0$ yields 
\[
0=B_{14}(\overline{x},\overline{y})=(y_1(f_1+x_4))^2\Rightarrow y_1=0\text{ or }x_4=y_1.
\]
If $y_1=0$, then $0=B_{33}(\overline{x},\overline{y})=f_1^2y_1^2$, so we can see that $y_2=y_1=0$ by assumption and also $y_4=0$ due to \eqref{y4b11}, so all $y_i$ equal $0$.

On the other hand, if $x_4=f_1$ and $y_1\ne 0$, then we have
\[
0=B_{22}(\overline{x},\overline{y})=f_1^2(y_2+x_2y_1)^2;
\]
therefore we get $y_2=x_2y_1$ and
\[
 0=B_{24}(\overline{x},\overline{y})=f_1^2y_1^2(x_2+1)^2,
\]
so that $x_2=1$, which then implies $0=B_{44}(\overline{x},\overline{y})=f_1^2\beta y_1^2$, contradicting our assumptions.

At this point we return to the other possible case in \eqref{y1y2b12}, namely the case $y_1=y_2$. It leads to
\[
 0=B_{14}(\overline{x},\overline{y})=y_2^2(x_4+f_1+f_5x_3),
\]
i.e. $y_2=0$ or $x_4=f_1+f_5x_3$. But $y_1=y_2$ and \eqref{y4b11} already imply that in the former case all $y_i$ vanish, whereas in the latter case we can hence assume $y_2\ne0$ and $1+x_2+x_3=0$ from
\[
0=B_{33}(\overline{x},\overline{y})=f_1^2y_2^2(1+x_2+x_3)^2=0.
\]
The final step is then to look at $B_{44}(\overline{x},\overline{y})$, which is equal to $f_1^2\beta y_2^2$ and thus gives the desired contradiction.

The only remaining case is $x_1=1=y_3$. The first helpful observation is
\[
0=B_{11}(\overline{x},\overline{y})=y_4+y_1x_4+f_5x_2+f_5x_3y_2,
\]
hence we must have 
\begin{equation}\label{y4rep}
 y_4=y_1x_4+f_5x_2+f_5x_3y_2.
\end{equation}
Using this consequence we obtain 
\[
 x_2=1+x_3(y_1+y_2)
 \]
 from
 \[
0=B_{12}(\overline{x},\overline{y})=f_5(x_2+1+x_3(y_1+y_2))^2.
\]
Thus we deduce that 
\begin{equation}\label{b34long}
0=B_{34}(\overline{x},\overline{y})=y_1x_3(x_4+f_5x_3+f_1y_1+f_1y_2)^2,
\end{equation}
i.e. $y_1=0$ or $x_3=0$ or $x_4=f_5x_3+f_1y_1+f_1y_2$. We will handle these cases separately.

Let us first suppose that $y_1=0$, in which case 
\begin{equation}\label{b14}
0=B_{14}(\overline{x},\overline{y})=f_5^2x_3^2(y_2+1)^2
\end{equation}
and thus $x_3=0$ or $y_2=1$.

In case $y_2=1$, we consider $K(\overline{y})=(y_4+f_5)^2$, so that $y_4=f_5$ and moreover
\[
 0=B_{33}(\overline{x},\overline{y})=(x_4+f_1+f_5x_3)^2
\]
implies $0=B_{22}(\overline{x},\overline{y})=\beta$ which cannot happen by assumption.

But if $x_3=0$ and $y_2\ne 1$, then we observe $0=B_{33}(\overline{x},\overline{y})=(x_4+f_1y_2)^2$, hence $0=B_{23}(\overline{x},\overline{y})=f_1(1+y_2)^2$ gives us a contradiction.

We proceed by assuming that $x_3=0\ne y_1$ in \eqref{b34long}; here we observe $0=B_{14}(\overline{x},\overline{y})=y_1^2(x_4+f_1)^2$, hence $x_4=f_1$. We then have $0=B_{33}(\overline{x},\overline{y})=f_1^2(1+y_1+y_2)^2$, so that we can deduce $y_1+y_2+1=0$ and thus $0=B_{22}(\overline{x},\overline{y})=\beta$, a contradiction.

The upshot of this is that in order to finish the proof of the lemma we can assume we are in the case $x_1=1=y_3$, $x_3y_1\ne0$ and $x_4=f_5x_3+f_1y_1+f_1y_2$ (see \eqref{b34long}). We can see immediately that
\[
0=B_{23}(\overline{x},\overline{y})=f_1(1+y_1+y_2)^2(1+x_3y_1)^2.
\]
Upon noticing 
\[
1+y_1+y_2=0\Rightarrow 0=B_{24}(\overline{x},\overline{y})=\beta x_3y_1
\]
we may thus assume that $x_3y_1=1$ and $y_1+y_2\ne 1$. 

We have
\[
0=B_{14}(\overline{x},\overline{y})=(1+y_1+y_2)^2(f_5x_3+f_1y_1)^2,
\]
resulting in $f_5x_3=f_1y_1$. This relation allows us to obtain
\[
 x_4=f_1y_2
\]
from \eqref{b34long} and hence $y_4=f_1y_1y_2$ from \eqref{y4rep}. We also have $f_5=f_5x_3y_1=f_1y_1^2$. Now we make these substitutions in 
$K(\overline{y})$ and find
\[
 0=K(\overline{y})=y_1^2(f_1^2y_2^4+f_1y_2^2+f_3+f_3^2)
\]
so $f_1^2y_2^4=f_1y_2^2+f_3+f_3^2$. But if we plug this into $B_{24}(\overline{x},\overline{y})$ we see that
\[
 0=B_{24}(\overline{x},\overline{y})=f_1(y_1+y_2+1)^2,
\]
contradicting the assumption $y_1+y_2+1\ne 0$.

This finally completes the proof of the lemma.



\begin{thebibliography}{999}
{\frenchspacing


  \bibitem{CasselsFlynn}
    {\sc J.W.S. Cassels}, {\sc E.V. Flynn}, \emph{Prolegomena to a middlebrow
    arithmetic of curves of genus $2$}, (Cambridge University Press, Cambridge, 1996).

  \bibitem{SD}
  \sc S.Duquesne}, {\em Traces of the group law on the Kummer
		surface of a curve of genus 2 in characteristic 2}, Preprint (2007).
 
  \bibitem{SD2}
   {\sc S. Duquesne}, `Montgomery scalar multiplication for genus 2 curves',
\emph{ANTS VI}, Burlington, VT, 2004, ed. D.A. Buell, Lecture Notes in Comput. Sci., vol. 3076 (Springer, Berlin-Heidelberg, 2004) 153--168.

	\bibitem{SD3}
  {\sc S. Duquesne}, `Montgomery ladder for all genus 2 curves in characteristic 2', \emph{WAIFI 2008}, Lecture Notes in Comput. Sci., vol. 5130 (Springer, Berlin, 2008) 174--188.

  \bibitem{Flynna}
   {\sc E.V. Flynn}, `The jacobian and formal group of a curve of genus 2 over an arbitrary ground field', 
    \emph{Math. Proc. Camb. Phil. Soc.} 107,  (1990) 425--441.

  \bibitem{Flynnb}
   {\sc E.V. Flynn}, `The group law on the jacobian of a curve of genus 2', 
    \emph{J. reine angew. Math.} 439, (1995), 45--69.

  \bibitem{FlynnSmart}
   {\sc E.V. Flynn}, {\sc N.P. Smart}, `Canonical heights on the jacobians
    of curves of genus $2$ and the infinite descent', 
    \emph{Acta Arith.} 79, (1997), 333--352.
  
  \bibitem{Gaudry1}
  {\sc P. Gaudry}, `Fast genus 2 arithmetic based on Theta functions',
   \emph{J. Math. Crypt.} 1, (2007), 243--265.
   
  \bibitem{GaudryLubicz}
   {\sc P. Gaudry}, {\sc D. Lubicz}, `The arithmetic of characteristic 2 Kummer
surfaces and of elliptic Kummer lines', \emph{Finite Fields Th. App.} 15, (2009) 246--260.    

   
 \bibitem{Hudson}
   {\sc R.W.H.T. Hudson}, \emph{Kummer's Quartic Surface} (University Press, Cambridge, 1905).

 \bibitem{Lang}
   {\sc S. Lang}, \emph{Introduction to Algebraic and Abelian Functions}, 2nd edition (Springer-Verlag, New York, 1982). 

 \bibitem{Magma}
  {\sc MAGMA} is described in {\sc W. Bosma}, {\sc J. Cannon} and {\sc C. Playoust}, 
  `The Magma algebra system I: The user language', 
  \emph{J. Symb. Comp.} 24, (1997), 235--265. (See also the Magma home page at \\{\tt http://magma.maths.usyd.edu.au/magma/}.)
  
  \bibitem{Maple}
   {\tt http://www.maplesoft.com/}.     
    
  \bibitem{StollH2}
   {\sc M.~Stoll}, `On the height constant for curves of genus two, II',
    \emph{Acta Arith.} 104, (2002), 165--182
\end{thebibliography}
\end{document}